\documentclass[a4paper]{article}
\newtheorem{theorem}{Theorem}
\newtheorem{lemma}{Lemma}
   
\newcommand{\EndOfProof}{\rule{0.5em}{1ex}}
\newcommand{\Strut}{\rule{0cm}{3ex}}

\title{On the Parametric Instability\\ Caused by Step Size Variation\\
in Runge-Kutta-Nystr\"om Methods\footnote{
Originally published as Mathematics Report~73, Tampere University of Technology
Department of Information Technology, isbn 9521500212,
issn 1236956x}
}
\author{Robert Pich\'e\\
   Tampere University of Technology\\
   Tampere, Finland\\
   {\tt robert.piche@tut.fi}}
\date{}
\begin{document}
\maketitle
\begin{abstract}
The parametric instability arising
when ordinary differential equations (ODEs)
are numerically integrated 
with Runge-Kutta-Nystr\"om (RKN) methods with varying step sizes is investigated.
Perturbation methods are used to 
quantify the critical step sizes associated with parametric instability.
It is shown that there is no parametric instability for
 linear constant coefficient ODEs  integrated with
RKN methods that are based on A-stable Runge-Kutta methods, 
because the solution is nonincreasing in some norm for all positive step sizes, 
constant or varying.
\end{abstract}
%
\section{Introduction}
In a recent paper~\cite{W98}, Wright showed how instability can arise when 
linear second-order constant-coefficient  ordinary differential equations (ODEs)
  are numerically integrated with varying step sizes, 
even when all steps are smaller than the critical step size
that ensures stability  in constant-step computations.
One of his examples was the central difference method
applied to the model problem
\begin{equation}   \label{linord2}
 \ddot{x}+x=0
\end{equation}
With constant step size $h$, 
the central difference method is known to be stable when $h<2$. 
Wright showed how it can become unstable when the step size
has small-amplitude oscillation of period 2 about the constant value $h=\sqrt{2}$.
The instability manifests itself as oscillation with growing amplitude.
Instability also arises if the step size
has period-3 oscillation about $h=1$, and
so on to period-6 oscillation about $h\approx 0.518$.
He conjectured that
instability could arise for  arbitrarily small step size.
Shortly after Wright's paper was published, Skeel
wrote a letter to the editor~\cite{S98} pointing out that  he had already proved
Wright's conjecture in a short note published in BIT~\cite{S93}.

In this work the instability phenomenon is investigated
by two approaches: first, stability analysis using contractivity concepts,
and second, parametric instability analysis using perturbation methods.
Here a brief description of these approaches is given.

In the usual stability analysis of a numerical ODE integration method,
the linear ODE system is transformed to diagonal form using a 
similarity transformation.
Because the transformation depends on the the step size,
this approach is only useful for methods used with constant step sizes.
Contractivity analysis is not based on the diagonalising transformation.
The idea is, given an ODE problem whose solution
is nonincreasing in some norm, to derive conditions under which
the numerical solution will also be nonincreasing.
Hairer et al.~\cite{HBL82} have used contractivity concepts
to derive stability results for a wide class of ODE integration methods
with constant step size.
The same approach is used here to derive stability results for
methods with varying step sizes.
It is shown that A-stable methods are stable also with variable
step size.

The second line of inquiry is prompted by Wright's
observation that
\begin{quote}
\ldots integrating a variable-stiffness system with constant time
steps resembles integrating a constant-stiffness system
with varying time steps.
\end{quote}
It is well known that oscillatory variation of
the stiffness parameter in variable-stiffness systems
can induce a kind of instability called {\em parametric resonance}~\cite{NM79}.
Wright's remark suggests that the
instability he observed is a kind of parametric resonance
of the ``discrete time system'' comprised by the numerical integration method
applied to the ODE.
Indeed, Wright's stability region diagrams 
resemble the classic Strutt diagrams used in analysis of parametric
resonance.

Although the literature on parametric resonance is vast,
there appears to be little on parametric resonance of difference equations.
Tanaka and Sato~\cite{TS93} have studied a second order difference equation
with periodic coefficients, a kind of discrete Mathieu equation.
They compute Strutt diagrams and show how to approximate the
stability region boundaries using perturbation methods.
This type of analysis is applied here to RKN integrators
with periodically varying step size.

\section{Runge-Kutta-Nystr\"om methods}
In this section a well-known general class of numerical integration methods
for second-order ODE initial value problems is presented.
Five specific methods from this class are singled out for further study.

The $s$-stage Runge-Kutta-Nystr\"om (RKN) method for advancing the
solution of the second-order velocity-independent ODE
\begin{equation}   \label{ord2}
    \ddot{x}=f(t,x)
\end{equation}
from $t=t_n$ to $t_{n+1}=t_n+h_n$ is
\begin{eqnarray*}
k_i&=&f(t_n+c_ih_n,x_n+c_i h_n\dot{x}_n+h_n^2\sum_{j=1}^s\bar{a}_{ij}k_j)
\;\;(1\leq i \leq s) \\
x_{n+1}&=&x_n+h_n\dot{x}_n+h_n^2\sum_{i=1}^s \bar{b}_ik_i\\
\dot{x}_{n+1}&=&\dot{x}_n+h_n\sum_{i=1}^s b_i k_i
\end{eqnarray*}
For example, the coefficients for
an explicit fourth-order RKN method~\cite[p.262]{HNW86} are
\begin{equation} \label{nystrom}
 \begin{array}{c|c}
    c_i&\bar{a}_{ij}  \\ \hline
     &\Strut \bar{b}_i  \\  
    & b_i
    \end{array}
=
 \begin{array}{c|ccc}
    0 \\
  1/2 & 1/8     \\
  1 & 0 & 1/2    \\ \hline
     & 1/6 & 1/3 & 0 \\ 
     & 1/6 & 4/6 & 1/6
  \end{array}
\end{equation}

A RKN method is equivalent to a Runge-Kutta method if there
exists a coefficient matrix $a_{ij}$ such that
\[
\bar{a}_{ij}=\sum_{k=1}^s a_{ik}a_{kj},\;
\bar{b}_i=\sum_{j=1}^s b_j a_{ji}
\]
For example, the A-stable third-order 2-stage SDIRK
method~\cite[p.203]{HNW86} corresponds to the RKN 
method~\footnote{
In~\cite{SFB90} it is shown that no
R-stable third-order diagonally implicit RKN method 
can have $-0.547 < c_1 < 1.213$,
but this is only for RKN methods
satisfying the simplifying condition
$c_j^2=2\sum_k a_{jk}$.}
\begin{equation}  \label{sdirk}
\begin{array}{c|cc}
    \alpha& \alpha^2 \\
    1-\alpha & 2\alpha-4\alpha^2 & \alpha^2    \\ \hline
     & (1-\alpha)/2 & \alpha/2 \\
     & 1/2 & 1/2 
    \end{array}
\end{equation}
where $\alpha=(3 +\sqrt{3})/6$.

The Newmark method for (\ref{ord2}) 
is a RKN method with two coefficients, $\beta$ and $\gamma$,
that appear in the RKN tableau as follows:
\[ 
\begin{array}{c|cc}
    0& 0 & 0   \\
    1&   (1-2\beta)/2 & \beta \\ \hline
     & (1-2\beta)/2 & \beta \\ 
     & 1-\gamma & \gamma 
    \end{array}
\]
The central difference method is given by the
Newmark parameters
\begin{equation}  \label{centdiff}
   \beta=0,\;\gamma=\frac{1}{2} \end{equation}
The Newmark method that was shown by Wright
to exhibit parametric instability is specified
by the parameter values
\begin{equation}   \label{newmark}
  \beta=\gamma=\frac{1}{2}  \end{equation}
The only choice of Newmark parameters that gives
a Runge-Kutta method (the trapezoid rule) is 
\begin{equation}  \label{trap}
   \beta=\frac{1}{4},\;\gamma=\frac{1}{2}
\end{equation}

Applying the RKN method to the 
linear constant coefficient model problem~(\ref{linord2}) 
and eliminating the stage variables gives 
\begin{equation}   \label{transition} 
y_{n+1} =R(h_n)y_n
 \end{equation}
where the state vector $y_n$ is defined as 
\[
  y_n:=\left[\begin{array}{c}x_{n}
   \\
  \dot{x}_{n}\end{array}\right]
\]
and the one-step
transition (amplification) matrix $R$ is given by
\[ R(h_n)=
 \left[\begin{array}{cc}
    1-h_n^2\bar{b}^T(I+h_n^2\bar{A})^{-1}e &
    h_n-h_n^3\bar{b}^T(I+h_n^2\bar{A})^{-1}c\\
    -h_n b^T(I+h_n^2\bar{A})^{-1}e &
    1-h_n^2b^T(I+h_n^2\bar{A})^{-1}c 
  \end{array}\right]
\]
with $e$ being a vector of ones.
In particular, the transition matrices for the RKN methods
presented earlier are as follows.
\begin{itemize}
\item The explicit fourth-order RKN method given in Eqn~(\ref{nystrom}):
\[
R(h)=\left [\begin {array}{cc} 
   1 & h \\
  -h & 1
\end {array}\right ]
+h^2\left [\begin {array}{cc} 
   -1/2+h^2/24 & -h/6 \\\noalign{\medskip}
  h/6-h^3/96 & -1/2+h^2/24
\end {array}\right ]
\]
\item The A-stable third-order 2-stage SDIRK of Eqn~(\ref{sdirk}):
\begin{eqnarray*}
\lefteqn{R(h)=\left [\begin {array}{cc} 
   1 & h \\
  -h & 1
\end {array}\right ]}
\\
&&+\frac{h^2}{(1+h^2\alpha^2)^2}
\left [\begin {array}{cc} 
-(4h^2\alpha^3-h^2\alpha^2+1)/2 & -(h^2\alpha^3-\alpha+1)\alpha h
   \\ \noalign{\medskip}
(h^2\alpha^3-\alpha+1)\alpha h & (4h^2\alpha^3-h^2\alpha^2+1)/2
\end {array}
\right ]
\end{eqnarray*}
\item The Newmark method:
\[
R(h)=\left [\begin {array}{cc} 
   1 & h \\
  -h & 1
\end {array}\right ]
+\frac{h^2}{(1+h^2\beta)}
\left [\begin {array}{cc} -1/2 & -\beta h \\ \noalign{\medskip}
\gamma h/2 & -\gamma \end {array}
\right ]
\]
with, as special cases, the central difference method of Eqn~(\ref{centdiff}):
\[
R(h)=\left [\begin {array}{cc} 
   1 & h \\
  -h & 1
\end {array}\right ]
+h^2
\left [\begin {array}{cc} -h/2 & 0 \\ \noalign{\medskip}
h/4 & -1/2 \end {array}
\right ]
\]
the Newmark method of Eqn~(\ref{newmark}):
\[
R(h)=\left [\begin {array}{cc} 
   1 & h \\
  -h & 1
\end {array}\right ]
+\frac{h^2}{(1+h^2/2)}
\left [\begin {array}{cc} -1/2 & -h/2 \\ \noalign{\medskip}
h/4 & -1/2 \end {array}
\right ]
\]
and the trapezoid method of Eqn~(\ref{trap}):
\[
R(h)=\left [\begin {array}{cc} 
   1 & h \\
  -h & 1
\end {array}\right ]
+\frac{h^2}{(1+h^2/4)}
\left [\begin {array}{cc} -1/2 & -h/4 \\ \noalign{\medskip}
h/4 & -1/2 \end {array}
\right ]
\]
\end{itemize}

\section{Contractivity Analysis \label{contractivity}}
%

%
%
%

The main result of this section is the following theorem,
which states that A-stable Runge-Kutta methods such as the trapezoid rule
are contractive in a norm that is independent of step size.
The same result of course holds for RKN methods
derived from A-stable Runge-Kutta methods.
Thus, these methods do not exhibit the
parametric instability caused by varying  step sizes.

\begin{theorem}
For any A-stable Runge-Kutta method applied to a
stable linear constant-coefficient ODE
there exists a
symmetric positive definite matrix $W$, independent of the
 step size,  such that the numerical
solution is contractive in the $W$-weighted euclidean norm,
that is, 
any two numerical solutions $y_n,\bar{y}_n$ satisfy
\[ \|y_{n+1}-\bar{y}_{n+1}\|_W \leq \|y_n-\bar{y}_n\|_W \]
\end{theorem}

Before looking at the proof of the theorem,
a couple of remarks. First, Wright's results for the
implicit R-stable Newmark method of Eqn~(\ref{newmark})
is a counterexample showing that
R-stability or implicitness are not
sufficient to preclude parametric resonance
with varying step size.

Secondly, in his letter to the editor, Skeel~\cite{S98}
mentions that the Newmark method of Eqn~(\ref{newmark})
\begin{quote}
{\ldots}does not satisfactorily deal with variable step size.
A method that does
is the implicit midpoint method, for which the amplification matrix
is an orthogonal matrix. Of course, it has the same
practical drawback of implicitness\ldots
\end{quote}
The implicit midpoint method has the same transition matrix
as the trapezoid rule of Eqn~(\ref{trap}). 
Because the transition matrix is orthogonal,
the method is contractive in the euclidean norm.
Thus, theorem~1 extends Skeel's remark to cover all
RKN methods that are equivalent to A-stable Runge-Kutta methods.
However, since A-stable methods are implicit,
this result is of limited use to those who want to use explicit methods.

The theorem is proved using three lemmas.
The first lemma states that 
any stable linear constant-coefficient system is 
contractive in some norm.
\begin{lemma}
If $A$ is Hurwitz then there exists  
a constant symmetric
positive definite matrix $W$
such that the solution of the ODE $\dot{y}=Ay$
is contractive in the $W$-weighted euclidean norm.
\end{lemma}
{\sc Proof.} 
It suffices to consider the homogeneous problem 
\[ \dot{y}=Ay \]
Since $A$ is Hurwitz, there exists
a symmetric positive definite matrix $W$
that satisfies the Lyapunov equation 
\[ A^TW+WA=-I  \]
Apply the linear transformation $z=W^{1/2}y$
to take the original ODE into 
\[\dot{z}=W^{1/2}AW^{-1/2}z\]
The solutions to the ODE are monotonically nonincreasing
in the euclidean norm, because 
\[ \frac{d}{dt} \| z \|^2 =-z^TW^{-1}z\leq 0 \]
Since $\|z\|=\|y\|_W$,
the stated result follows.~\EndOfProof
\vspace{2ex}

The $s$-stage Runge-Kutta method for advancing the
solution of the linear constant
coefficient ODE
\[ \dot{y}=Ay\]
from $t_n$ to $t_{n+1}=t_n+h$ is given by
\begin{eqnarray*}
g_i&=&y_n + h\sum_{j=1}^s a_{ij}Ag_j\;\;\;(1\leq i\leq s)\\
y_{n+1}&=&y_n+h\sum_{j=1}^s b_jAg_j
\end{eqnarray*}
Eliminating the stage variables $g_i$ results in
\[ y_{n+1}=R(hA)y_n \]
where $R$ is the stability function of the Runge-Kutta method.

A linear transformation $z=Vy$, where $V$ is a constant
nonsingular matrix, gives a new ODE
\[ \dot{z}=VAV^{-1}z \]
The next result tells us that the result of applying
the Runge-Kutta method to the transformed ODE 
is the same as transforming the numerical results 
from the original ODE.
\begin{lemma}
The following diagram commutes:
\[ \begin{array}{ccccc}
    &\mbox{\fbox{$y_n$}} & 
     \stackrel{V}{\longrightarrow} & 
    \mbox{\fbox{$z_n$}}  \\
    R(hA) & \downarrow &&\downarrow &R(hVAV^{-1}) \\
   & \mbox{\fbox{$y_{n+1}$}} & 
    \stackrel{V}{\longrightarrow} &
     \mbox{\fbox{$z_{n+1}$}}
\end{array}
\]
\end{lemma}
{\sc Proof.}
Substitution.~\EndOfProof
\vspace{2ex}

Finally, the following result from~\cite[Corollary 11.3]{HW91} is needed.
\begin{lemma} 
If an A-stable Runge-Kutta method
is applied to a linear constant-coefficient ODE 
that is contractive in the euclidean norm,  
then the numerical solution is contractive in the
euclidean norm.
\end{lemma}

Combining the three lemmas gives the theorem.
Lemma~1 gives the weight $W$ under which the  ODE is
contractive.
$V=W^{1/2}$ is the linear transformation matrix that
relates the original variables $y$ to the new variables $z$
that are contractive in the euclidean norm.
By lemma~3, numerical sequences $z_n$  computed with A-stable
Runge-Kutta methods inherit the contractivity of the ODE,
and by lemma~2, this also applies to the original ODE.

\section{Stability charts}

A necessary and sufficient
condition for stable integration of the model problem~(\ref{linord2}) 
with constant step size
is  given
by the Schur-Cohn inequalities
\begin{equation}   \label{SC}
| \mbox{trace}(R) | -1 \leq \mbox{det}(R) \leq 1
\end{equation} 
Thus, with constant step sizes the
 central difference method, Eqn~(\ref{centdiff}),
is stable for $0\leq h\leq 2$,
the RKN method of Eqn~(\ref{nystrom})
is stable for $0\leq h\leq 2(2+2^{1/3}-2^{2/3})^{1/2}\approx 2.59 $,
and the SDIRK method of Eqn~(\ref{sdirk}),
the trapezoid method (Eqn~(\ref{trap})), and 
the Newmark~ method of Eqn~(\ref{newmark}) 
are stable for all $h\geq 0 $, i.e.\ they are R-stable.

Now consider integration  of the model problem~(\ref{linord2}) with
a sequence of step sizes $h_n$
that is periodic with period $p$.
The composition of $p$  steps can be interpreted
as a one-step method integrating from $t=t_n$ to 
$t=t_{n+p}=t_n+h_n+\cdots+h_{n+p-1}$.
The application of one step of the composed method
gives
\[  y_{n+p} 
 =R([h_n,\ldots,h_{n+p-1}]) y_n
\]
where 
the transition matrix $R$ of the composed method is given by
the matrix product
\[ R([h_n,\ldots,h_{n+p-1}]):=
     R(h_{n+p-1})R(h_{n+p-2})\cdots R( h_{n+1})R( h_n )
\]
This composed transition matrix determines the stability
of the linear homogeneous difference equation~(\ref{transition}).
The stability condition is given
by the inequalities~(\ref{SC}) applied with the
composed transition matrix.

For the periodic step size variation given by
\begin{equation}   \label{hn}
   h_n=h+\epsilon\cos(2\pi n/p)
\end{equation}
a stability chart can be made
by numerically testing the stability criterion
for a large number of specific values of $(h,\epsilon)$. 
Figure~1 shows values giving an unstable composed 
Nystr\"om method of Equation~(\ref{nystrom})
when the step size period is $p=6$.
Also shown are the corresponding stability charts
for the central difference method and the Newmark
method of Equation(\ref{newmark}).

The instability regions for the
central difference and Newmark methods
form wedges with their points on the $h$-axis.
Stability charts for other periods $p$ are similar,
except that the number of wedges increases with $p$.
The figure resembles 
a Strutt diagram for the Mathieu equation~\cite{NM79}.

The instability regions for the Nystr\"om method
are more rounded and, for values of $h$ smaller than
the constant-step stability limit, do not reach the $h$-axis. 
This figure resembles Strutt diagrams for damped oscillators.

The stability charts 
for the trapezoid and SDIRK methods
are not shown, because they contain no instability points.
This is because they are A-stable Runge-Kutta methods,
which 
according to the theorem in
section~\ref{contractivity} remain stable with varying step size.

\section{Perturbation Analysis of Stability Regions}
The points on the stability charts
where the instability regions 
intersect the $h$ axis
will be called {\em critical}\/ step sizes.
In this section,
their values are identified using perturbation analysis.

Substituting the trial solution
\begin{equation}   \label{trialsolution}
 y_n=y^{(0)}_n+\epsilon y^{(1)}_n+\epsilon^2  y^{(2)}_n+\cdots
\end{equation}
and the step size (\ref{hn}) into (\ref{transition})
and equating powers of $\epsilon$ gives a sequence of 
coupled difference equations, the first three of which are
\begin{eqnarray}
y^{(0)}_{n+1}&=&R(h) y^{(0)}_n   \label{eps0}
   \\
y^{(1)}_{n+1}&=&R(h) y^{(1)}_n + \cos(2\pi n/p)R'(h)y^{(0)}_{n}
     \label{eps1} 
   \\
y^{(2)}_{n+1}&=&R(h) y^{(2)}_n + \cos(2\pi n/p)
   \left(R'(h)y^{(1)}_n + \frac{1}{2} R''(h)y^{(0)}\right)
      \nonumber
\end{eqnarray}
The first equation, Equation~(\ref{eps0}), is a homogeneous linear
constant coefficient difference equation, and is stable provided that
the constant component $h$ of the step size is
within the stability range for the constant-step method.

The remaining equations
are linear constant coefficient difference
equations with a harmonic forcing term.
If the eigenvalues of $R$ are 
smaller than one in magnitude
(i.e.\ the method has algorithmic damping),
then these equations will all have bounded solutions.
For such methods, the ``wedges''
in the stability chart do not reach the $h$ axis,
and the method is stable 
when the 
amplitude $\epsilon$ of step size oscillation
is sufficiently small.
This is confirmed by the stability chart (Figure~1) for
the Nystr\"om method of Eqn~(\ref{nystrom}),
which has algorithmic damping.

Resonance instability may arise in integration methods
without algorithmic damping
if the forcing frequency matches the natural
frequency.
It therefore suffices to restrict attention to 
methods
whose transition matrix $R(h)$ has
eigenvalues of unit modulus.
These eigenvalues can be written as
\begin{equation}  \label{eigs}
   \mbox{eig}(R(h)) = \exp( \pm i \omega(h)) 
\end{equation}
where $\omega$ is a real function of $h$.
 For example, the central difference method
has
\[ \omega(h) = \arctan(h\sqrt{1-h^2/4},1-h^2/2)\;\;(0\leq h\leq 2)\]
while the Newmark method of Eqn~(\ref{newmark}) has
\[
\omega(h) = \arctan(h\sqrt{1+h^2})\;\;(0\leq h)   
\]

For methods satisfying (\ref{eigs}), the solution of (\ref{eps0})
is a linear combination of trigonometric functions of the form
\begin{equation}  \label{y0}
 y^{(0)}_n= Y_1 \cos(n\omega) +
Y_2\sin(n\omega)
\end{equation}
Substituting (\ref{y0})  
into (\ref{eps1}) 
gives a forcing function for $y^{(1)}_n$ of the form
\[ R'(h)Y_1 \cos(n\omega)\cos(2\pi n/p) +
   R'(h)Y_2 \sin(n\omega)\cos(2\pi n/p)
\]
Resonance occurs when
\begin{equation}   \label{criticalh}
   \omega(h) = \frac{\pi}{p}
\end{equation}
or
\[
  \omega(h)=\pi-\frac{\pi}{p}
\]
Values of $h$ satisfying these equations are critical
 step sizes. 

The smallest critical step size for period-$p$ 
 step size oscillation
is denoted $h^0$, and is 
given by the solution of Equation~(\ref{criticalh}).
The minimum critical step sizes for various $p$ are listed
in Tables~\ref{centdiffcriticalh} and
\ref{newmarkcriticalh}.
These values agree with those found by Wright~\cite{W98}
by other means.

\begin{table}
\centering
\begin{tabular}{c|r@{$\,\approx\,$}l | r@{$\,\approx\,$}l}
  $p$ & 
 \multicolumn{2}{c|}{$h^0$} &
 \multicolumn{2}{c}{$h^1$} \\ \hline
2& $\sqrt{2}$ & 1.4142 & \multicolumn{2}{c}{$\pm 1/2$} \\
3 & \multicolumn{2}{c|}{1} & \multicolumn{2}{c}{$\pm 1/8$} \\
4 & $\sqrt{2-\sqrt{2}}$ & 0.7654 &  $\pm(2-\sqrt{2})/8$ & $\pm$0.0732 \\
5 &$(\sqrt{5}-1)/2$ & 0.6180 & $\pm(3-\sqrt{5})/16$ & $\pm$0.0477 \\
6 & $(\sqrt{6}-\sqrt{2})/2$ & 0.5176& 
   $\pm(4-\sqrt{12})/16$ & $\pm$0.0335
\end{tabular}
\caption{\label{centdiffcriticalh}
Smallest critical  step size $h^0$ 
and wedge width $h^1$ for period-$p$
oscillation, central difference method}
\end{table}

\begin{table}
\centering
\begin{tabular}{c|r@{$\,\approx\,$}l | r@{$\,\approx\,$}l}
  $p$ & 
 \multicolumn{2}{c|}{$h^0$} &
 \multicolumn{2}{c}{$h^1$} \\ \hline
2& \multicolumn{2}{c}{$\infty$} & \multicolumn{2}{c}{--} \\
3 & $\sqrt{2}$ & 1.4142 & \multicolumn{2}{c}{$\pm$1/4} \\
4 &  $\sqrt{-2+2\sqrt{2}}$ & 0.9102& $\pm(-1+\sqrt{2})/4$ & 0.1036 \\
5 & $\sqrt{-4+2\sqrt{5}}$ & 0.6871 & $\pm(-2+\sqrt{5})/4$ & 0.0590 \\
6 & $\frac{1}{3}\sqrt{-18+12\sqrt{3}}$ & 0.5562 & $\pm(-3+2\sqrt{3})/12$ & 0.0387
\end{tabular}
\caption{\label{newmarkcriticalh}
Smallest critical  step size $h^0$ 
and wedge width $h^1$ for period-$p$
oscillation, Newmark method}
\end{table}

Any RKN method of order at least 1 has
\[ R(h)= \left[\begin{array}{cc}
    1 & h \\ -h &1
  \end{array}\right]+O(h^2)
\]
and so
\[   \mbox{eig}(R(h)) = \pm ih + O(h^2) \]
The smallest critical
step size for large $p$ is then given by
\[ \min h_{\mbox{crit}} \approx \frac{\pi}{p}   \]
Thus, the
critical step size can be arbitrarily small.
This agrees with Skeel's result~\cite{S93} for the central difference method.

The method of strained parameters is now used to
analyse the shape of the instability region ``wedge''
near the $h$-axis.
The average step size parameter is assumed to have the form
\[
h=h^0+h^1\epsilon + O(\epsilon^2)
\]
where $h^0$ is a critical step size.
The step size then varies according to
\[ h_n=h^0+h^1\epsilon+\epsilon\cos(2\pi n/p) +O(\epsilon^2) \]
Substituting this into the composed transition matrix
gives
\[  R([h_n,\ldots,h_{n+p-1}]) = R_0+\epsilon R_1 + O(\epsilon^2) \]
where
\[ R_0:=R(h^0)^3\]
and
\[ R_1:=(\cos(2\pi 0/p)+h^1)R(h^0)^{p-1}R'(h^0)
   +\ldots+(\cos(2\pi (p-1)/p)+h^1)R'(h^0)R(h^0)^{p-1}\]
The boundaries of the instability regions
are found by solving the stability conditions for the composed
transition matrix. 

For example, for the central difference method
with step size period $p=3$ we have $h^0=1$ and
\[
R([h_2,h_1,h_0])=-I+\frac{\epsilon}{16}
\left[\begin{array}{cc} 
   6 & -64h^1-4 \\ 48h^1-3 & -6  
   \end{array}\right] + O(\epsilon^2)
\]
Neglecting terms in $\epsilon^2$, 
the stability condition~(\ref{SC}) is
\[ 1 \leq 1-\frac{3}{16}{\epsilon}^{2}+12\,{\epsilon}^{2}{(h^1)}^{2}
  \leq 1  \]
which is satisfied when
\[   h^1=\pm \frac{1}{8} 
\]
The instability region boundary is therefore approximated by
\begin{equation}    \label{wedge}
  h = 1 \pm \frac{1}{8}\epsilon + O(\epsilon^2)
\end{equation}
Figure~2 shows instability points in the neighbourhood
of the critical step size for the central difference
method for parametric oscillation of period $p=3$,
and the instability region boundaries predicted by~(\ref{wedge}).
The approximation appears to be correct.
Tables~\ref{centdiffcriticalh} and
\ref{newmarkcriticalh} give instability region ``wedge''
widths for the central difference method and the Newmark method.

\section{Conclusions}
The following results were presented:
\begin{itemize}
\item A-stable Runge-Kutta methods are stable when applied with
varying time step to stable linear time invariant ODE systems.
\item For numerical ODE integration methods that have numerical damping,
there will be no parametric resonance if the step size oscillation amplitude
is sufficiently small.
\item The critical average step size about which small oscillation may cause
parametric resonance in the numerical integration of a simple oscillator
can be found using perturbation methods.
\item The perturbation analysis indicates that there can be arbitrarily small
critical step sizes.
\item The method of strained parameters can be used
to approximate the shape of the
instability region in Strutt diagrams for small amplitude
step size oscillation.
\end{itemize}
Further investigation is needed for the following questions:
\begin{itemize}
\item Is A-stability necessary to rule out parametric resonance?
\item What about integration of nonlinear ODEs? 
\item How can a ``safe'' step size oscillation amplitude 
be estimated for ODE integration methods that have numerical damping? 
\end{itemize}


\end{document}